\documentclass[11pt]{amsart}
\usepackage{amsrefs}

\newcommand{\eps}{\epsilon}
\newcommand{\pa}{\partial}

\newfont{\fnt}{cmr10 scaled 550}
\renewcommand{\eps}{\varepsilon} 

\newtheorem{theorem}{Theorem}

\newtheorem{conj}{Conjecture}
\newtheorem{lemma}{Lemma}
\newtheorem{cor}{Corollary}
\newtheorem{prop}{Proposition}

\newtheorem{definition}{Definition}

\theoremstyle{remark}
\newtheorem{remark}{Remark}

\font\strange=msbm10

\renewcommand{\epsilon}{\varepsilon}

\newcommand{\C}{{{\mathchoice  {\hbox{$\textstyle{\text{\strange C}}$}}
{\hbox{$\textstyle{\text{\strange C}}$}}
{\hbox{$\scriptstyle{\text{\strange C}}$}}
{\hbox{$\scriptscriptstyle{\text{\strange C}}$}}}}}

\renewcommand{\Sigma}{\varSigma}
\newcommand{\R}{{{\mathchoice  {\hbox{$\textstyle{\text{\strange R}}$}}
{\hbox{$\textstyle{\text{\strange R}}$}}
{\hbox{$\scriptstyle  N\kern-0.3em  R$}}  
{\hbox{$\scriptscriptstyle  R\kern-0.2em  R$}}}}}

\newcommand{\Z}{{{\mathchoice  {\hbox{$\textstyle{\text{\strange Z}}$}}
{\hbox{$\textstyle{\text{\strange Z}}$}}
{\hbox{$\scriptstyle  Z\kern-0.3em  Z$}}
{\hbox{$\scriptscriptstyle  Z\kern-0.2em  Z$}}}}}

\newcommand{\N}{{{\mathchoice  {\hbox{$\textstyle{\text{\strange N}}$}}
{\hbox{$\textstyle{\text{\strange N}}$}}
{\hbox{$\scriptstyle  N\kern-0.3em  N$}}
{\hbox{$\scriptscriptstyle  N\kern-0.2em  N$}}}}}

\renewcommand{\phi}{\varphi}
\usepackage{amsmath,amsthm,amscd}

\begin{document}
\title{Normal scalar curvature conjecture and its applications}
\date{November 14, 2007}

 \author{Zhiqin Lu}

\address{Department of
Mathematics, University of California,
Irvine, Irvine, CA 92697}


\email[Zhiqin Lu]{zlu@math.uci.edu}

\thanks{The
 author is partially supported by  NSF award DMS-09-04653.}

\begin{abstract}
In this paper, we proved the Normal Scalar Curvature Conjecture and the B\"ottcher-Wenzel Conjecture.
We developed a new Bochner formula
and it becomes useful with the first  conjecture we proved.
Using the results, we  established some new  pinching theorems for minimal submanifolds in spheres.
\end{abstract}

\maketitle 

\tableofcontents 
\pagestyle{myheadings}

\newcommand{\ka}{K\"ahler }
\newcommand{\ii}{\sqrt{-1}}
\section{Introduction}
Let $M^n$ be an  $n$-dimensional manifold isometrically immersed into the space form $N^{n+m}(c)$ of constant sectional curvature $c$. The normalized scalar curvature $\rho$  is defined as follows:
\begin{equation}
\rho=\frac{2}{n(n-1)}\sum_{1=i<j}^n R(e_i, e_j,e_j,e_i),\\
\end{equation}
where $\{e_1,\cdots, e_n\}$  is a local  orthonormal frame of the tangent bundle, and $R$ is the curvature tensor for the tangent  bundle.

The (normalized) scalar curvature of the normal bundle is defined as:
\[
\rho^\perp=\frac{1}{n(n-1)}|R^\perp|,
\]
where $R^\perp$ is the curvature tensor of the normal bundle.
More precisely, let $\xi_1,\cdots,\xi_m$ be a local orthonormal frame of the normal bundle. Then
\begin{equation}
\rho^{\perp}=\frac{2}{n(n-1)}\left(\sum_{1=i<j}^n\sum_{1=r<s}^m\langle
R^\perp(e_i,e_j)\xi_r,\xi_s\rangle^2\right)^{\frac 12}.
\end{equation}
Unlike the normalized scalar curvature, $\rho^\perp$ is always nonnegative.

In the study of submanifold theory, De Smet, Dillen, Verstraelen, and Vrancken~\cite{ddvv} proposed  the following {\sl Normal Scalar Curvature Conjecture}\footnote{Also known as the DDVV conjecture.}: 

\begin{conj}\label{cc1} Let $h$ be the second fundamental form, and 
let $H=\frac 1n \,{\rm trace}\, h$ be the mean curvature tensor.  Then 
\begin{equation}\label{0p}
\rho+\rho^\perp\leq |H|^2+c.
\end{equation}
\end{conj}

In the first part of this paper, we proved the above conjecture.

\smallskip

 Let $x\in M$ be a fixed point and let $(h_{ij}^r)$ ($i,j=1,\cdots,n$ and $r=1,\cdots,m$) be the entries  of (the traceless part of) the second fundamental form under the  orthonormal frames of both the tangent bundle and the normal bundle. Then by
 ~\cites{dfv,bogd}, Conjecture~\ref{cc1} can be formulated as an inequality with respect to the entries $(h_{ij}^r)$ as follows:
  \begin{align}\label{1a}
\begin{split}
&\sum_{r=1}^m\sum_{1=i<j}^n(h_{ii}^r-h_{jj}^r)^2+2n\sum_{r=1}^m\sum_{1=i<j}^n(h_{ij}^r)^2\\
&\geq 2n\left(\sum_{1=r<s}^m\sum_{1=i<j}^n\left(\sum_{k=1}^n(h_{ik}^rh_{jk}^s-h_{ik}^sh_{jk}^r)\right)^2\right)^{\frac 12}.
\end{split}
\end{align}

Let $A$ be an $n\times n$ matrix. Let
\[
||A||=\sqrt{\sum_{i,j=1}^n a_{ij}^2}
\]
be its Hilbert-Schmidt norm,
where $(a_{ij})$ are the entries of $A$.
Let
\[
[A,B]=AB-BA
\]
be the commutator of two $n\times n$ matrices. 
Suppose that $A_1,A_2,\cdots,A_m$ are $n\times n$ symmetric real matrices. 
Then  inequality~\eqref{1a}, in terms of matrix notations, can be formulated as:
\begin{conj}\label{conj2}
 For $n, m\geq 2$, we have
\begin{equation}\label{conj}
(\sum_{r=1}^m||A_r||^2)^2\geq 2(\sum_{r<s}||[A_r,A_s]||^2).
\end{equation}
\end{conj}

Fixing $n,m$, we call the above inequality Conjecture $P(n,m)$. Conjecture ~\ref{cc1} is equivalent to Conjecture $P(n,m)$ for any positive numbers $n,m$.

 A weaker version of  Conjecture~\ref{cc1},
$
 \rho\leq |H|^2+c,
 $
 was proved in ~\cite{ch1}. An alternate proof is in~\cite{bogd2}.

The following special cases of Conjecture~\ref{cc1} were known. $P(n,2)$ was proved in~\cite{chern-d-k}; $P(2,m)$ was proved in~\cite{ddvv}; $P(3,m)$ was proved in~\cite{cl}; and $P(n,3)$ was proved in~\cite{lu21}. In~\cite{dfv}, a weaker version of $P(n,m)$ was proved by using an algebraic inequality in~\cite{lianmin} (see also~\cite{xusenlin}) . In the same paper, $P(n,m)$ was proved under the addition assumption that the submanifold is either {\it Lagrangian $H$-umbilical}, or {\it ultra-minimal} in $\mathbb C^4$. In~\cite{dfv-2},  Conjecture~\ref{cc1} was studied for  invariant submanifolds of K\"ahler and Sasakian manifolds and in particular, it was proved for complex submanifolds of  $\C^n$.  Finally, in~\cite{ge-tang}, an independent (and different) proof of Conjecture~\ref{cc1} was given.

The proofs of Conjectures ~\ref{cc1} and~\ref{bw-conj} were also given in~\cite{lu22}, the previous version of this paper.

It should be pointed out that the classification of the submanifolds when the equality in~\eqref{0p} holds is a very difficult problem. 
An easy and special case was done in~\cite{cl}. More systematically, the problem was treated in the recent preprint of 
Dajczer and Tojeiro~\cite{dajczer-tojeiro}.

In the second part of this paper, we used the method in proving Conjecture~\ref{cc1}
to sharpen the pinching theorems of Simons type~\cite{simons-j}.  The inequality of Simons was improved by many people (for an incomplete list, ~\cites{chern-d-k,lianmin,xusenlin}). 
By their works,  it is well known that for an $n$-dimensional manifold $M$ minimally immersed into $S^{n+m}$, we have: 1). If $m=1$ and $0<||\sigma||^2\leq n$, then $||\sigma||^2=n$ and $M=M_{r,n-r}$; 2). If $m>1$, and $0<||\sigma||^2\leq\frac 23 n$, then $||\sigma||^2=\frac 23n$ and $M$ has to be the Veronese surface\footnote{For the definition of $M_{r,n-r}$ and the Veronese surface, see \S~\ref{se}.}.

In the past, people studied the Laplacian of the norm of the second fundamental form. However, more accurate results will be obtained by studying the Laplacian of the  norm of the second fundamental form on {\it each} normal direction.
We   established new Simons-type  formula~\eqref{s-71} for the above idea. The key linear algebraic inequality ~\eqref{key-1} in proving Conjecture ~\ref{cc1} is just the right tool to make the formula  useful.

 We got  a new pinching theorem (Theorem~\ref{bhu}). The theorem unified  and sharpened the previous pinching theorems, and may become the starting point of the
gap theorem of Peng-Terng~\cite{peng-terng} type in high codimensions (see Conjecture~\ref{pt}).

In the last part of this paper, we  proved the conjecture of B\"ottcher and Wenzel~\cite{bw}. The conjecture was from the theory of random matrices and is purely linear algebraic in nature.

\begin{conj}[B\"ottcher-Wenzel Conjecture]\label{bw-conj}
Let $X,Y$ be two $n\times n$ matrices. Then
\begin{equation}\label{conj-3-1}
||[X,Y]||^2\leq 2||X||^2||Y||^2.
\end{equation}
\end{conj}

In ~\cite{bw}, the following weaker version of the conjecture was proved.
\[
||[X,Y]||^2\leq 3||X||^2||Y||^2.
\]
 
 \smallskip 
 
 The proof of the conjecture was posted in~\cite{lu22}. Shortly after that, there are  two independent and different proofs of the conjecture in~\cites{vj-1,bw-2}.
 \\
 
 To get an idea of the proofs of Conjectures~\ref{cc1} and ~\ref{bw-conj}, we first observe the following theorem~\cite{chern-d-k}*{Lemma 1}, which proves $P(n,2)$:
 
 \begin{theorem}\label{poi}
  Let $A,B$ be $n\times n$ symmetric matrices. Then
 \begin{equation}\label{chern-1}
 ||[A,B]||^2\leq 2||A||^2||B||^2.
 \end{equation}
 \end{theorem}

 {\bf Proof.} Without loss of generality,  we  assume that $A$ is a diagonal matrix. Let
 \[
 A=\begin{pmatrix}\lambda_1&\\&\ddots&\\&&\lambda_n\end{pmatrix}.
 \]
 Let $B=(b_{ij})$. Then
 \begin{equation}\label{wer}
 ||[A,B]||^2=2\sum_{i,j}(\lambda_i-\lambda_j)^2b_{ij}^2.
 \end{equation}
 The theorem follows from the fact that
 \[
 (\lambda_i-\lambda_j)^2\leq 2(\lambda_i^2+\lambda_j^2)\leq  2\sum_k\lambda_k^2.
 \]
 
 \qed
 
 The above result can be viewed as a baby version of both Conjectures ~\ref{cc1} and ~\ref{bw-conj}. In fact, from the above inequality, we get
 \[
 \sum_{k=2}^s||[A_1,A_k]||^2\leq 2||A_1||^2(\sum_{k=2}^s||A_k||^2)
 \]
 for any symmetric matrices $A_1,\cdots,A_s$. The key step in proving Conjecture~\ref{cc1} is a refinement of the above inequality into  the following version:
 \[
 \sum_{k=2}^s||[A_1,A_k]||^2\leq ||A_1||^2(\sum_{k=2}^s||A_k||^2+\underset{2\leq k\leq s}{\rm Max}\, ||A_k||^2).
 \]
 The  inequality is new even when $s=2$.
  See Remark~\ref{remk1} for more details.
  
  In addition to  the above,
 a trick in proving Conjecture~\ref{bw-conj} is as follows:
 if we let
 \[
 {\rm ad}\,(A)\, B=[A,B],
 \]
and if $A$ is diagonalized. Then the eigenvalues of the operator ${\rm ad}(A)$, acting on the space $n\times n$ matrices, have multiplicity at least $2$. In fact, let $0\neq B=(b_{ij})$ be a symmetric matrix such that
\[
{\rm ad}\,(A)\,B=\lambda B.
\]
Define $B'=(b_{ij}')$, where $b_{ij}'=b_{ij}$ for $i>j$ and $b_{ij}'=-b_{ij}$ for $i<j$. Then $B'$ is also an eigenvector of the same eigenvalue.

If $A$ is not a symmetric matrix. We found that  $B'= [A^T,B^T]$ serves the same purpose. This is one of the crucial step in
the proof.

Finally, we can generalize ~\eqref{conj-3-1} into the following infinite dimensional version\footnote{We thank Timur Oikhberg for the help in the infinite dimensional setting. There is  an infinite dimensional version also for  the Normal Scalar Curvature Conjecture in terms of linear algebraic inequalities.}, which can be proved by  operator approximation by matrices.

\begin{theorem} Let $H$ be a separable Hilbert space and let $A,B$ be linear operators with finite Hilbert-Schmidt norms. Then we have
\[
||[A,B]||^2\leq 2||A||^2\cdot||B||^2.
\]
\end{theorem}

\qed

 {\bf Acknowledgement.}  The author thanks Professor Chuu-Lian Terng for the helps in many aspects during the preparation of this paper. He also thanks Bogdan Suceav{\u{a}} for his bringing the Normal Scalar Curvature Conjecture into his attention, without which the paper is not possible.

\section{Invariance}

Let $A_1,\cdots,A_m$ be  $n\times n$ symmetric matrices.
Let $G=O(n)\times O(m)$. Then $G$ acts on matrices $(A_1,\cdots,A_m)$ in the following natural way: let $(p,q)\in G$, where $p,q$ are $n\times n$ and $m\times m$ orthogonal matrices, respectively. Let $q=\{q_{ij}\}$. Then
\[
(p,I)\cdot (A_1,\cdots,A_m)=(pA_1p^{-1},\cdots, pA_mp^{-1}),
\]
and 
\[
(I,q)\cdot (A_1,\cdots,A_m)=(\sum_{j=1}^m q_{1j}A_j,\cdots,\sum_{j=1}^mq_{mj}A_{j}).
\]

It is easy to verify the following:

\begin{prop}\label{prop1}
Conjecture $P(n,m)$ is $G$ invariant. That is, in order to prove inequality ~\eqref{conj}
for $(A_1,\cdots,A_m)$, we just need to prove the inequality for any $\gamma\cdot (A_1,\cdots,A_m)$ where $\gamma\in G$. Moreover, the expressions of both sides of ~\eqref{conj} are $G$ invariant.
\end{prop}

\qed

\begin{cor}\label{cor11}
We can prove Conjecture~\ref{conj2} under the following additional assumptions on the matrices:
\begin{enumerate}
\item $A_1$ is diagonal;
\item $\langle A_\alpha,A_\beta\rangle=0$ if $\alpha\neq \beta$;
\item $||A_1||\geq\cdots\geq ||A_m||$.
\end{enumerate}
\end{cor}

\qed

Note that under the above assumptions, $A_k=0$ if $k>\frac 12n(n+1)$.

\section{Proof of the normal scalar curvature conjecture}

In this section, we prove Conjecture~\ref{conj2}. We first establish some lemmas which are themselves interesting.

\begin{lemma} \label{lem1-1}
Suppose $\eta_1,\cdots,\eta_n$ are real numbers and 
\[
\eta_1+\cdots+\eta_n=0,\quad
\eta_1^2+\cdots+\eta_n^2=1.
\]
Let $r_{ij}\geq 0$ be nonnegative numbers for $i<j$. Then we have
\begin{equation}\label{use}
\sum_{i<j}(\eta_i-\eta_j)^2 r_{ij}\leq\sum_{i<j}r_{ij}+{\rm Max} (r_{ij}).
\end{equation}
If $\eta_1\geq\cdots\geq\eta_n$, and $r_{ij}$ are not simultaneously  zero, then the equality in~\eqref{use} holds in one of the following three cases:
\begin{enumerate}
\item $r_{ij}=0$ unless $(i,j)=(1,n)$, $(\eta_1,\cdots,\eta_n)=(1/\sqrt 2,0,\cdots,0, -1/\sqrt 2)$;
\item $r_{ij}=0$ if $2\leq i<j$, $r_{12}=\cdots=r_{1n}\neq 0$, and\newline  $(\eta_1,\cdots,\eta_n)=(\sqrt{(n-1)/n},-1/\sqrt{n(n-1)},\cdots,-1/\sqrt{n(n-1)})$;
\item $r_{ij}=0$ if $i<j<n$, $r_{1n}=\cdots=r_{(n-1)n}\neq 0$, and\newline
 $(\eta_1,\cdots,\eta_n)=(1/\sqrt{n(n-1)},\cdots,1/\sqrt{n(n-1)}, -\sqrt{(n-1)/n})$;
\end{enumerate}
\end{lemma}

{\bf Proof.}
We assume that
$
\eta_1\geq\cdots\geq\eta_n.
$
If $\eta_1-\eta_n\leq 1$ or $n=2$, then ~\eqref{use} is trivial. So we assume $n>2$,  and
\[
\eta_1-\eta_n>1.
\]
We observe that $\eta_i-\eta_j< 1$ for $2\leq i<j\leq n-1$. Otherwise, we could have
\[
1\geq\eta_1^2+\eta_n^2+\eta_i^2+\eta_j^2\geq\frac 12
((\eta_1-\eta_n)^2+(\eta_i-\eta_j)^2)>1,
\]
which is a contradiction. 

Using the same reason,
if  $\eta_1-\eta_{n-1}>1$,  then we have $\eta_2-\eta_{n}\leq 1$; and if $\eta_2-\eta_{n}>1$, then we have $\eta_1-\eta_{n-1}\leq 1$.
Replacing $\eta_1,\cdots,\eta_n$ by $-\eta_n,\cdots,-\eta_1$ if necessary, we can always assume that $\eta_2-\eta_{n}\leq 1$.
Thus $\eta_i-\eta_j\leq 1$ if $2\leq i<j$, and ~\eqref{use} is implied by the following inequality
\begin{equation}\label{use2}
\sum_{1<j}(\eta_1-\eta_j)^2 r_{1j}\leq\sum_{1<j} r_{1j}+\underset{1<j}{\rm Max}\, (r_{1j}).
\end{equation}
Let $s_j=r_{1j}$ for $j=2,\cdots,n$. 
Then the above inequality becomes
\begin{equation}\label{use3}
\sum_{1<j}(\eta_1-\eta_j)^2 s_{j}\leq\sum_{1<j} s_{j}+\underset{1<j}{\rm Max}\, (s_{j}).
\end{equation}

In order to prove ~\eqref{use3}, we define the matrix $P$ as follows
\[
P=\begin{pmatrix}
\underset{{1<j}}{\sum} s_j& -s_{2}&\cdots &-s_{n}\\
-s_{2}& s_{2}\\
\vdots&&\ddots
\\
-s_{n}&&&s_{n}
\end{pmatrix}.
\]

We 
claim that the maximum eigenvalue of $P$ is no more than $r=\sum_j s_j+{\rm Max}\, (s_j)$. To see this, we compute the determinant of the matrix
\[
\begin{pmatrix}
y-\underset{{1<j}}{\sum} s_j& s_{2}&\cdots &s_{n}\\
s_{2}& y-s_{2}\\
\vdots&&\ddots
\\
s_{n}&&&y-s_{n}
\end{pmatrix}
\]
for $y>r$.
Using the Cramer's rule,  the answer is
\[
(y-s_2)\cdots (y-s_n)\left( y-\sum_{1<j}s_j-\sum_{1<j}\frac{s_j^2}{y-s_j}\right).
\]
We have $y-s_k>\sum_{j=2}^n s_j$ for any $1<k\leq n$. Thus the above expression is greater than
\[
(y-s_2)\cdots (y-s_n)(y-\sum_{1<j}s_j-(\sum_{1<j}s_j)^{-1}\sum_{1<j} s_j^2)>0.
\]

Let $\eta=(\eta_1,\cdots,\eta_n)^T$, we then have
\[
\sum_{1<j}(\eta_1-\eta_j)^2 s_{j}=
\eta^TP\eta\leq r=\sum_{1<j}^n s_j+\underset{1<j}{\rm Max}\, (s_j).
\]

We assume that $n>2$.
If the equality in~\eqref{use} holds, then we must have $\eta_1-\eta_n>1$. Otherwise
\[
\sum_{i<j}(\eta_i-\eta_j)^2 r_{ij}\leq\sum_{i<j}r_{ij}
\]
and all $\{r_{ij}\}$'s have to be zero. Since $\eta_1-\eta_n>1$, then $\eta_i-\eta_j<1$ for $2\leq i<j<n$. Thus $r_{ij}=0$ for $2\leq i<j<n$. Moreover, the equality of ~\eqref{use3} must hold. From the proof of~\eqref{use3}, we conclude that either at most one of $s_j$'s can be nonzero, or all $s_j$'s are the same. Translating this fact to $r_{ij}$, we conclude that if $r_{1n}\neq 0$, then either $r_{1j}=0$ for $j<n$, or $r_{12}=\cdots =r_{1n}\neq 0$. 
In the first case, there are two possibilities: either $r_{1n}=\cdots=r_{(n-1)n}\neq 0$, or
$r_{2n}=\cdots=r_{(n-1)n}=0$.
Putting the information together, we conclude that only in the three cases in the lemma the equality holds. This completes the proof.

\qed

\begin{lemma}\label{lem2}
Let $A$ be an $n\times n$  diagonal  matrix of norm $1$. Let $A_2,\cdots, A_m$ be symmetric matrices such that
\begin{enumerate}
\item $\langle A_\alpha,A_\beta\rangle=0$ if  $\alpha\neq \beta$;
\item $||A_2||\geq\cdots\geq||A_m||$.
\end{enumerate}
Then we have
\begin{equation}\label{key-1}
\sum_{\alpha=2}^m||[A,A_\alpha]||^2\leq \sum_{\alpha=2}^m ||A_\alpha||^2+||A_2||^2.
\end{equation}
The  equality in~\eqref{key-1} holds if and only if, after an orthonormal base change and up to a sign, we have
\begin{enumerate}
\item
$A_3=\cdots =A_m=0$, and
\begin{equation}\label{a116}
A_1=\begin{pmatrix}
\frac{1}{\sqrt 2} &0\\0&-\frac{1}{\sqrt 2}\\
&&0\\
&&&\ddots\\
&&&&0
\end{pmatrix},\quad
A_2=c\begin{pmatrix}
0&\frac{1}{\sqrt 2}\\\frac{1}{\sqrt 2}&0\\
&&0\\
&&&\ddots\\
&&&&0
\end{pmatrix},
\end{equation}
where $c$ is any constant, or
\item For two real numbers $\lambda=1/\sqrt{n(n-1)}$ and $\mu$, we have
\begin{equation}\label{a112}
A_1=\lambda\begin{pmatrix}
n-1\\
&-1\\
&&\ddots\\
&&&-1
\end{pmatrix},
\end{equation}
and
$A_\alpha$ is $\mu$ times the matrix whose only nonzero entries are $1$ at the $(1,\alpha)$ and $(\alpha,1)$ places, where $\alpha=2,\cdots,n$.
\end{enumerate}
\end{lemma}

{\bf Proof.} We assume that each $A_\alpha$ is not zero.
Let $A_\alpha=((a_\alpha)_{ij})$, where  $(a_\alpha)_{ij}$ are  the entries for $\alpha=2,\cdots,m$. Let
\[
\delta=\underset{i\neq j}{\rm Max}\,\sum_{\alpha=2}^m(a_\alpha)^2_{ij}.
\]
Let 
\[
A=\begin{pmatrix}
\eta_1\\
&\ddots\\
&&\eta_n
\end{pmatrix}.
\]
Then by the previous lemma, we have
\begin{equation}\label{ab-1}
\sum_{\alpha=2}^m||[A,A_\alpha]||^2\leq \sum_{\alpha=2}^m||A_\alpha||^2+2\delta.
\end{equation}
Thus it remains to prove that
\begin{equation}\label{fun}
2\delta\leq ||A_2||^2.
\end{equation}
To see this, we identify each $A_\alpha$ with the (column)   vector $\vec A_\alpha$ in $\mathbb R^{\frac 12 n(n+1)}$
as follows:
\[
A_\alpha\mapsto (a_{12},\cdots,a_{1n},a_{23},\cdots,a_{2n},\cdots,a_{(n-1)n},
\frac{1}{\sqrt 2}a_{11},\cdots,\frac{1}{\sqrt 2}a_{nn})^T.
\]
 Let $\mu_\alpha$ be the norm of the vector $A_\alpha$. Then we have
\begin{equation}\label{cd-1}
\mu_\alpha^2=\frac 12 ||A_\alpha||^2
\end{equation}
for $\alpha=2,\cdots,m$. Extending the set of vectors $\{\vec A_\alpha/\mu_\alpha\}_{2\leq\alpha\leq m}$ into an orthonormal basis of $\mathbb R^{\frac 12n(n+1)}$
\[
\vec A_2/\mu_2, \cdots, \vec A_m/\mu_m,\vec A_{m+1},\cdots, \vec A_{\frac 12n(n+1)+1},
\]
we get an orthogonal matrix. Apparently, each row vector of the matrix is a unit vector. Thus we have
\[
\sum_{\alpha=2}^m (\mu_\alpha)^{-2}(a_\alpha)^2_{ij}\leq 1
\]
for fixed $i<j$.
Since $\mu_2\geq\cdots\geq\mu_m$, we get
\[
\sum_{\alpha=2}^m (a_\alpha)^2_{ij}\leq\mu_2^2\leq \frac 12||A_2||^2.
\]
This proves ~\eqref{fun}. Finally, when equality holds, 
according to Lemma~\ref{lem1-1}, there are three cases. 
The first case corresponds to the first case in Lemma~\ref{lem2}.
The second and the third cases in Lemma~\ref{lem1-1} are equivalent by the permutation $(\eta_1,\cdots,\eta_n)\to
(-\eta_n,\cdots,-\eta_1)$.  Translating to the notations in Lemma~\ref{lem2}, $A_1$ is in the form of~\eqref{a112}. Moreover, we have 
\[
(a_\alpha)_{ij}=0
\]
for $\alpha=2,\cdots,n$, and $1<i<j$. Since $A_1$ is invariant under the similar transformation $A_1\mapsto QA_1Q^T$, where $Q$ is of the form
\[
\begin{pmatrix}
1\\
&Q_1
\end{pmatrix},
\]
and $Q_1$ is an orthogonal matrix. Up to an orthonormal base change and up to a sign, $A_2,\cdots,A_n$ can be represented as in the second case 
of Lemma~\ref{lem2}. This completes the proof.

\qed

\begin{remark}\label{remk1}
Let $A$ be a diagonal matrix of unit norm and let $B$ be a symmetric matrix. Let $||B||_\infty={\rm Max}\, (|b_{ij}|)$, where $(b_{ij})$ are the entries of $B$. By~\eqref{ab-1}, we get
\[
||[A,B]||^2\leq ||B||^2+2||B||_\infty^2.
\]
Although not directly used in this paper, this is a sharper estimate. 
Note that Theorem~\ref{poi} is a much weaker version of the above inequality. This shows that, even in the case of $m=2$, Lemma~\ref{lem2} is a refinement of previous results.
\end{remark}

{\bf Proof of Conjecture~\ref{conj2}.} 
Let $a>0$ be the largest positive real number such that
\[
(\sum_{\alpha=1}^m||A_\alpha||^2)^2\geq 2a(\sum_{\alpha<\beta}||[A_\alpha,A_\beta]||^2).
\]
Since $a$ is maximum, 
by Corollary~\ref{cor11}, we can find matrices $A_1,\cdots,A_m$ such that
\begin{equation}\label{conj-1}
(\sum_{\alpha=1}^m||A_\alpha||^2)^2= 2a(\sum_{\alpha<\beta}||[A_\alpha,A_\beta]||^2)
\end{equation}
with the following additional properties:
\begin{enumerate}
\item $A_1$ is diagonal;
\item $\langle A_\alpha,A_\beta\rangle=0$ if $\alpha\neq \beta$;
\item $0\neq ||A_1||\geq ||A_2||\geq\cdots \geq ||A_m||$.
\end{enumerate}

We let $t^2=||A_1||^2$ and let $A'=A_1/|t|$. Then ~\eqref{conj-1} becomes a quadratic expression in terms of $t^2$:
\begin{align*}&
t^4-2t^2(a\sum_{1<\alpha}||[A',A_\alpha]||^2-\sum_{1<\alpha}||A_\alpha||^2)
+(\sum_{\alpha=2}^m||A_\alpha||^2)^2\\&\qquad- 2a(\sum_{1<\alpha<\beta}||[A_\alpha,A_\beta]||^2)=0.
\end{align*}
Since the left-hand side of the above is nonnegative for all $t^2$, we have
\[
a\sum_{1<\alpha}||[A',A_\alpha]||^2-\sum_{1<\alpha}||A_\alpha||^2>0,
\]
and 
\[
||A_1||^2=a\sum_{1<\alpha}||[A',A_\alpha]||^2-\sum_{1<\alpha}||A_\alpha||^2.
\]
By  Lemma~\ref{lem2}, we have
\[
\sum_{1<\alpha}||[A',A_\alpha]||^2\leq \sum_{\alpha=2}^m||A_\alpha||^2+||A_2||^2\leq \sum_{\alpha=1}^m||A_\alpha||^2,
\]
which proves that $a\geq 1$.

\qed

\section{The optimal inequality}
Let $A_1,\cdots,A_m$ be $n\times n$ symmetric matrices. Assume that $||A_1||=\cdots =||A_m||=1$. Let
\[
\sigma_{ij}=||[A_i,A_j]||^2
\]
for $i,j=1,\cdots,m$. 

From~\eqref{conj}, we get the following result:
\begin{prop}\label{pp-1}
Let $x_1,\cdots,x_m\geq 0$ be nonnegative real numbers. Then we have
\[
\sum_{i,j=1}^m \sigma_{ij}x_i x_j\leq(\sum_{i=1}^m x_i)^2.
\]
\end{prop}

We make the following definition:

\begin{definition}
A symmetric $m\times m$ matrix $P=(p_{ij})$ is called pseudo-positive, if for any nonnegative real numbers $x_1,\cdots,x_m\geq 0$, 
\[
\sum_{i,j=1}^m p_{ij} x_i x_j\geq 0.
\]
\end{definition}

A symmetric matrix has property $K$, if 
for any negative eigenvalue of the matrix, the components of the corresponding eigenvector are neither all nonpositive nor all nonnegative. Using the Lagrange's multiplier's method,
we can characterize the pseudo-positiveness as follows:

\begin{prop}\label{pp-2}
$A$ is a pseudo-positive matrix if and only if any principal submatrix of $A$ has property $K$.
\end{prop}

Using the above notations, we can reformulate Proposition~\ref{pp-1} as follows:

\begin{prop}\label{pp-3}
Let $\Sigma=(\sigma_{ij})$ and let $S$ be the $m\times m$ matrix whose entries are all $1$. Then
$S-\Sigma$ is a pseudo-positive matrix.
\end{prop}

\qed

The main result of this section is to show that the Normal Scalar Curvature Conjecture implies the following result in
~\cite{lianmin}*{p.585, Eq. (5)}~\footnote{The proof of ~\cite{xusenlin} is more geometric.}:
 
\begin{theorem}\label{lili}
 Let $x_1,\cdots,x_m\geq 0$ be nonnegative real numbers. Then
\begin{equation}\label{s-6}
\sum_{i,j=1}^m \sigma_{ij} x_i x_j\leq\frac 32 
(\sum_{i=1}^m x_i)^2-\sum_{i=1}^m x_i^2.
\end{equation}
\end{theorem}

{\bf Proof.} We use math induction.  Assume that the inequality~\eqref{s-6} is true for $m-1$. Let $x_1$ be the largest number among $x_1,\cdots,x_m$. Then we can rewrite equation~\eqref{s-6} as follows:
\begin{equation}\label{s-7-1}
\frac 12 x_1^2+x_1(3\sum_{j=2}^m x_j-2\sum_{j=2}^m \sigma_{1j}x_j)
+\frac 32(\sum_{j=2}^m x_j)^2-\sum_{j=2}^m x_j^2-\sum_{i,j=2}^m \sigma_{ij} x_i x_j\geq 0.
\end{equation}

If 
\[
3\sum_{j=2}^m x_j-2\sum_{j=2}^m \sigma_{1j}x_j\geq 0,
\]
then~\eqref{s-7-1} is true by the inductive assumption. Otherwise, the left-hand side of ~\eqref{s-7-1} attains its minimal when
\[
x_1=2\sum_{j=2}^m \sigma_{1j} x_j-3\sum_{j=2}^m x_j.
\]
Since $\sigma_{1j}\leq 2$ by Theorem~\ref{poi}, we have
\[
x_1\leq 4\sum_{j=2}^m x_j-3\sum_{j=2}^m x_j=\sum_{j=2}^m x_j.
\]
Since $x_1$ is the largest number among nonnegative numbers $x_1,\cdots,x_m$, we have
\[
\sum_{j=1}^m x_j^2\leq \frac 12 (\sum_{j=1}^m x_j)^2.
\]
By Proposition~\ref{pp-3}, we have
\[
\sum_{i,j=1}^m\sigma_{ij}x_i x_j\leq(\sum_{j=1}^m x_j)^2\leq \frac 32 
(\sum_{i=1}^m x_i)^2-\sum_{i=1}^m x_i^2.
\]
\qed

\section{Pinching theorems}\label{se}

Let $M$ be an $n$-dimensional compact minimal submanifold in the unit sphere $S^{n+m}$ of dimension $n+m$. 
Following~\cite{chern-d-k}, we make the following convention on the range of indices:
\begin{align*}
&1\leq A,B,C,\cdots\leq n+m;\quad 1\leq i,j,k,\cdots\leq n;\\
&n+1\leq \alpha,\beta,\gamma,\cdots\leq n+m.
\end{align*}

Let $\omega_{1},\cdots,\omega_{n+m}$ be an orthonormal frame of the cotangent bundle of $S^{n+m}$.
Then
we
 have
\begin{align}\label{s-1}
\begin{split}
&d\omega_A=-\omega_{AB}\wedge\omega_B,
\\
&d\omega_{AB}=-\omega_{AC}\wedge\omega_{CB}+\frac 12
K_{ABCD}\omega_C\wedge\omega_D,
\end{split}
\end{align}
where $\omega_{AB}$ are the connection forms and $K_{ABCD}$ is the curvature tensor of the sphere
\begin{equation}\label{s-1-1}
K_{ABCD}=\delta_{AC}\delta_{BD}-\delta_{AD}\delta_{BC}.
\end{equation}

Let $\omega_1,\cdots,\omega_n$ be an orthonormal frame of $TM$ and let
$\omega_{n+1},\cdots,\omega_{n+m}$ be an orthonormal frame of $T^\perp M$. Then we have
\begin{align}\label{s-2}
\begin{split}
&d\omega_i=-\omega_{ij}\wedge\omega_j,\\
&d\omega_{ij}=-\omega_{ik}\wedge\omega_{kj}+\frac 12 R_{ijkl}\omega_k\wedge\omega_l,
\end{split}
\end{align}
where $R_{ijkl}$ is the curvature tensor of $M$. We have the similar equations for the normal bundle:
\begin{equation}\label{s-3}
d\omega_{\alpha\beta}=-\omega_{\alpha\gamma}\wedge\omega_{\gamma\beta}+\frac 12
R_{\alpha\beta kl}\omega_k\wedge\omega_l,
\end{equation}
where $R_{\alpha\beta kl}$ is the curvature tensor of the normal bundle.

Comparing ~\eqref{s-1},~\eqref{s-2}, we have
\begin{equation}\label{s-4}
R_{ijkl}=K_{ijkl}+h_{ik}^\alpha h_{jl}^\alpha-h_{il}^\alpha h_{jk}^\alpha,
\end{equation}
where
\[
\omega_{\alpha i}=h^\alpha_{ij}\omega_j.
\]

Comparing ~\eqref{s-1} and~\eqref{s-3}, we have
\begin{equation}\label{s-5}
R_{\alpha\beta kl}=K_{\alpha\beta kl}+h_{ik}^\alpha h_{il}^\beta-h_{il}^\alpha h_{ik}^\beta.
\end{equation}

By~\eqref{s-1-1} from~\eqref{s-5},~\eqref{s-6}, we have
\begin{align*}
& R_{ijkl}=\delta_{ik}\delta_{jl}-\delta_{il}\delta_{jk}+h_{ik}^\alpha h_{jl}^\alpha-h_{il}^\alpha h_{jk}^\alpha,\\
&R_{\alpha\beta kl}=h_{ik}^\alpha h_{il}^\beta-h_{il}^\alpha h_{ik}^\beta.
\end{align*}

Define the covariant derivative of $h_{ij}^\alpha$ by
\begin{equation}\label{okj}
h_{ijk}^\alpha\omega_k=dh_{ij}^\alpha-h_{il}^\alpha\omega_{lj}-h_{lj}^\alpha\omega_{li}+h_{ij}^\beta\omega_{\alpha\beta}.
\end{equation}
Then we have
\begin{equation}\label{s-6-1}
h_{ijk}^\alpha=h_{ikj}^\alpha.
\end{equation}

Define the second covariant derivative of $h_{ij}^\alpha$ by
\[
h_{ijkl}^\alpha\omega_l=dh_{ijk}^\alpha-h_{ljk}^\alpha\omega_{li}-h_{ilk}^\alpha\omega_{lj}
-h_{ijl}^\alpha\omega_{lk}+h_{ijk}^\beta\omega_{\alpha\beta}.
\]
Then
\[
h_{ijkl}^\alpha-h_{ijlk}^\alpha=h_{ip}^\alpha R_{pjkl}+h_{pj}^\alpha R_{pikl}-h_{ij}^\beta R_{\alpha\beta kl}.
\]

Thus we have
\begin{equation}\label{psrq}
h^\alpha_{kijk}=h^\alpha_{kikj}+h^\alpha_{kp}R_{pijk}+h^\alpha_{pi}R_{pkjk}
-h^\beta_{ki}R_{\alpha\beta jk}.
\end{equation}

Define $\Delta h_{ij}^\alpha=h^\alpha_{ijkk}$. Then by~\eqref{s-6-1},~\eqref{psrq}, and the minimality of $M$, we have
\begin{align*}
&
\Delta h_{ij}^\alpha=h^\alpha_{kp}(\delta_{pj}\delta_{ik}-\delta_{pk}\delta_{ij}
+h^\beta_{pj}h^\beta_{ik}-h^\beta_{pk}h^\beta_{ij})\\
&+h^\alpha_{pi}((n-1)\delta_{pj}-h^\beta_{pk}h^\beta_{jk})-h^\beta_{ki}(h^\alpha_{pj}h^\beta_{pk}-h^\alpha_{pk} h^\beta_{pj})\\&
=nh_{ij}^\alpha+2 h^\alpha_{kp}h^\beta_{pj} h^\beta_{ik}-h_{kp}^\alpha h^\beta_{pk} h^\beta_{ij}
-h_{pi}^\alpha h_{pk}^\beta h_{jk}^\beta-h_{ki}^\beta h_{pj}^\alpha h_{pk}^\beta.
\end{align*}

Let $A^\alpha$ be the matrix of $h_{ij}^\alpha$. Then in terms of matrix notations, we have
\begin{equation}\label{s-7}
\Delta
A^\alpha=nA^\alpha-\langle A^\alpha, A^\beta\rangle
A^\beta-[A^\beta, [A^\beta, A^\alpha]].
\end{equation}

Before stating the theorems, we make the following definitions (from~\cite{chern-d-k}). For $1\leq r\leq n$, the submanifold $M_{r,n-r}$ is defined as
\[
M_{r,n-r}=S^r\left(\sqrt{\frac{r}{n}}\right)\times S^{n-r}\left(\sqrt{\frac{n-r}{n}}\right),
\]
which is immersed in $S^{n+1}$ in a natural way.  Since $S^{n+1}$ is a totally geodesic submanifold of  $S^{n+m}$, $M_{r,n-r}$ is regarded as a minimal submanifold of $S^{n+m}$ as well. The Veronese surface is defined as follows: let $(x,y,z)$ be the natural coordinate system in $\R^3$ and $(u^1,u^2,u^3,u^4,u^5)$ the natural coordinate system in $\R^5$. We consider the mapping defined by 
\begin{align*}
&
u^1=\frac{1}{\sqrt 3} yz,\,\, u^2=\frac{1}{\sqrt 3} zx, \,\,u^3=\frac{1}{\sqrt 3} xy, \,\,u^4=\frac{1}{2\sqrt 3}(x^2-y^2), \\
&u^5=\frac 16 (x^2+y^2-2z^2).
\end{align*}
This defines an isometric immersion of $S^2(\sqrt 3)$ into $S^4$. Since $S^4$ is totally geodesic in $S^{2+m}$,  the Veronese surface is a minimal surface of $S^{2+m}$.

Let $||\sigma||^2$ be the square of the length of the second fundamental form. Through the works of Simons~\cite{simons-j}, Chern-do Carmo-Kobayashi~\cite{chern-d-k}, Yau~\cite{yau-pinching}, Shen~\cite{shen-yibing}, and Wu-Song~\cite{wu-song}, and finally by Li-Li~\cite{lianmin}, Chen-Xu~\cite{xusenlin}, we get  the following  pinching theorem:

\begin{theorem}\label{thm303}
Let $M$ be an $n$-dimensional compact minimal submanifold in $S^{n+m}$, $m\geq 2$. If $||\sigma||^2\leq\frac 23n$ everywhere on $M$, then $M$ is either a totally geodesic submanifold or a Veronese surface in $S^{2+m}$.\footnote{There is a misprint in~\cite{lianmin}*{Theorem 3}. In fact, for any immersion $M\to S^{4}\to S^{2+m}$ for $m\geq 2$, $||\sigma||^2=4/3$.}
\end{theorem}

The proof is based on  the following Simons-type formula which can easily be derived from~\eqref{s-7}:

\begin{equation}\label{simo}
\frac 12\Delta||\sigma||^2=\sum_{i,j,k,\alpha}(h^\alpha_{ijk})^2+n||\sigma||^2-\sum_{\alpha,\beta}||[A_\alpha,A_\beta]||^2-\sum_{\alpha,\beta}|\langle A_\alpha,A_\beta\rangle |^2.
\end{equation}

Using Theorem~\ref{lili} and the maximal principle, we get $h_{ijk}^\alpha\equiv 0$, and $||\sigma||^2\equiv \frac 23 n$. Using ~\cite{chern-d-k}*{page 70}, we conclude that $M$ has to be either totally geodesic, or a Veronese surface.

The codimensional  $1$ case was studied in~\cite{chern-d-k}:

\begin{theorem} Let $M$ be a minimal hypersurface in $S^{n+1}$ such that
$0\leq||\sigma||^2\leq n$. Then $M$ is either totally geodesic, or one of $M_{r,n-r}$.
\end{theorem}

In this section, we sharpen the above results. Before stating the theorem, we make the following definition:

\begin{definition} The fundamental matrix $S$ of $M$ is an $m\times m$ matrix-valued function defined as $S=(a_{\alpha\beta})$, where
\[
a_{\alpha\beta}=\langle A^\alpha, A^\beta\rangle.
\]
We let $\lambda_1\geq\cdots\geq\lambda_m$ be the set of eigenvalues of the matrix. In particular, $\lambda_1$ is the largest eigenvalue and $\lambda_2$ is the second largest eigenvalue of the matrix $S$.
\end{definition}

Using the above notation, $||\sigma||^2$ is the trace of the fundamental matrix: $||\sigma||^2=\lambda_1+\cdots+\lambda_n$.
We have the following:

\begin{theorem}\label{bhu}
Let
\[
0\leq ||\sigma||^2+\lambda_2\leq n.
\]
Then $M$ is totally geodesic, or is one of $M_{r,n-r}$ ($1\leq r\leq n$) in $S^{n+m}$, or is a Veronese surface in $S^{2+m}$.
\end{theorem}

\begin{remark} Since
\[
\lambda_2\leq\frac 12||\sigma||^2,
\]
The theorem generalizes the above two theorems. 
\end{remark}

{\bf Proof.} For each integer $p\geq 2$, we define the smooth function~\footnote{At one point, $f_p=\sum\lambda_i^p$. However, it is in general not possible to find a smooth local frame such that the fundamental matrix is diagonalized on an open set. This is one of the technical difficulty of the theorem.}
\[
f_p={\rm Tr}\,(S^p).
\]
Let $x\in M$ be a fixed point. Let $(x_1,\cdots,x_n)$ be the local coordinates of $x$. We assume that at $x$, the fundamental matrix is diagonalized. A straightforward computation using~\eqref{s-7} gives that, at $x\in M$,
\begin{align}\label{s-71}
\begin{split}
&\frac{1}{2p}\Delta f_p=\frac 12\sum_{s+t=p-2}\sum_{k,\alpha,\beta}\lambda_\alpha^s\lambda_\beta^t\left(D_{\frac{\pa}{\pa x_k}} a_{\alpha\beta}\right)^2+\sum_\alpha\left(\lambda_\alpha^{p-1}\sum_{i,j,k}(h_{ijk}^\alpha)^2\right)\\
&+nf_p-f_{p+1}-\sum_{\beta\neq\alpha}||[A^\beta, A^\alpha]||^2||A^\alpha||^{p-1},
\end{split}
\end{align}
where
\[
D_{\frac{\pa}{\pa x_k}} a_{\alpha\beta}=\frac{\pa a_{\alpha\beta}}{\pa x_k}
+\omega_{\alpha\gamma}\left(\frac{\pa}{\pa x_k}\right) a_{\gamma\beta}
+\omega_{\beta\gamma}\left(\frac{\pa}{\pa x_k}\right) a_{\alpha\gamma}
\]
is the covariant derivative.

 We assume that at $x$,
\[
\lambda_1=\cdots=\lambda_r>\lambda_{r+1}\geq\cdots\geq\lambda_m.
\]
Then we have
\begin{align*}
&\frac 1p\Delta f_p\geq(p-1)\sum_{k,\alpha}\lambda_\alpha^{p-2}\left(D_{\frac{\pa}{\pa x_k}} a_{\alpha\alpha}\right)^2
+2\sum_\alpha\left(\lambda_\alpha^{p-1}\sum_{i,j,k}(h_{ijk}^\alpha)^2\right)
 \\&+2(nf_p-f_{p+1}-\sum_{\beta\neq\alpha}||[A^\beta, A^\alpha]||^2\lambda_\alpha^{p-1}).
\end{align*}
Using Lemma~\ref{lem2} and the above inequality, we get\footnote{If $r=m$, we define $A^{r+1}=0$.}
\begin{align}\label{nop}
\begin{split}
&\frac 1p\Delta f_p\geq(p-1)\sum_{k,\alpha}\lambda_\alpha^{p-2}\left(D_{\frac{\pa}{\pa x_k}} a_{\alpha\alpha}\right)^2
+2\sum_\alpha\left(\lambda_\alpha^{p-1}\sum_{i,j,k}(h_{ijk}^\alpha)^2\right)\\
&\quad +2( r||A^1||^{2p}(n-||A^1||^2-\sum_{\alpha=2}^m||A^\alpha||^2-\lambda_2))-6mn\lambda_{r+1}^{p}.
\end{split}
\end{align}

We have
\begin{equation}\label{fgh}
|\nabla f_p|^2=p^2\sum_{k}\left(\sum_\alpha\lambda_\alpha^{p-1}D_{\frac{\pa}{\pa x_k}} a_{\alpha\alpha}\right)^2.
\end{equation}

Using the Cauchy inequality, we get
\begin{equation}\label{nop-1}
|\nabla f_p|^2\leq p^2 f_p \sum_{k,\alpha}\lambda_\alpha^{p-2}\left(D_{\frac{\pa}{\pa x_k}} a_{\alpha\alpha}\right)^2.
\end{equation}

Let $g_p=(f_p)^{\frac 1p}$. Then at $f_p\neq 0$, using ~\eqref{nop} and~\eqref{nop-1}, we have
\begin{align}\label{kkk}
\begin{split}
&\Delta g_p=\frac 1p f_p^{\frac 1p  -1}\Delta f_p+\frac 1p (\frac 1p-1) f_p^{\frac 1p -2}|\nabla f_p|^2\\
&\geq  2f_p^{\frac 1p-1}\sum_\alpha\left(\lambda_\alpha^{p-1}\sum_{i,j,k}(h_{ijk}^\alpha)^2\right)\\&
 \quad+2f_p^{\frac 1p-1}( r||A^1||^{2p}(n-||A^1||^2-\sum_{\alpha=2}^m||A^\alpha||^2-\lambda_2)-3mn\lambda_{r+1}^{p}).
\end{split}
\end{align}

By~\eqref{fgh}, we have
\[
|\nabla g_p|\leq C||\sigma||
\]
for some constant $C$.
Thus we have
\[
\int_M\Delta g_p=0.
\]
Using this fact, from~\eqref{kkk}, we get
\begin{align*}
&\qquad \int_Mf_p^{\frac 1p-1} \sum_\alpha\left(\lambda_\alpha^{p-1}\sum_{i,j,k}(h_{ijk}^\alpha)^2\right)\\
&+\int_Mf_p^{\frac 1p-1} ( ||A^1||^{2p}(n-||A^1||^2-\sum_{\alpha=2}^m||A^\alpha||^2-\lambda_2)-3mn\lambda_{r+1}^{p})\leq 0.
\end{align*}
Since $\lambda_{r+1}^p/f_p\to 0$ almost everywhere when $p\to \infty$, from the above inequality, we get
\[
\int_M\sum_{i,j,k}\sum_{\alpha\leq r}(h^\alpha_{ijk})^2+||A^1||^2(n-||A^1||^2-\sum_{\alpha=2}^m||A^\alpha||^2-\lambda_2)
\leq 0.
\]
 Thus
 $\lambda_2+||\sigma||^2\equiv n$, $h^\alpha_{ijk}=0$ for $\alpha\leq r$, and we have
\[
\sum_{\alpha=2}^m||[A^1,A^\alpha]||^2\equiv ||A^1||^2(\sum_{\alpha=2}^m ||A_\alpha||^2+||A_2||^2).
\]
By Lemma~\ref{lem2}, there are four cases. 

{\bf Case 1.}  All $A^i$ are zero, then $M$ is totally geodesic. 

{\bf Case 2.} $A^2=\cdots =A^m=0$. In this case,  $||\sigma||^2=n$. Using~\eqref{simo}, we shall get
\[
h_{ijk}^\alpha=0
\]
for any $i,j,k=1,\cdots,n$ and $\alpha=n+1,\cdots,n+m$. 
Now we can use the techniques similar to those in~\cite{chern-d-k}.
With the suitable choice of local frame, we can assume that $h_{ij}^{n+1}=0$ if $i\neq j$. For any $i,j$ such that $h^{n+1}_{ii}\neq h^{n+1}_{jj}$, by~\eqref{okj}, we have
\[
0=d h_{ij}^{n+1}=(h^{n+1}_{ii}-h^{n+1}_{jj})\omega_{ij},
\]
and thus $\omega_{ij}=0$. From the structure equations, we get
\[
\frac 12R_{ijkl}\omega_k\wedge\omega_l=d\omega_{ij}+\omega_{ik}\wedge\omega_{kj}=0.
\]
Using  formula~\eqref{s-4}  for $R_{ijkl}$, we get
\[
(h^{n+1}_{ii} h^{n+1}_{jj}+1)\omega_i\wedge\omega_j=0.
\]
In particular, $\{h_{ii}^{n+1}\}$ can take at most two different values $\lambda_1,\lambda_2$ such that
$\lambda_1\lambda_2+1=0$. 

Let $r$ be the number of $\lambda_1$'s. Then from $r\lambda_1+(n-r)\lambda_2=0$ and $r\lambda_1^2+(n-r)\lambda^2=n$, we have $\lambda_1=\sqrt{\frac{n-r}{r}}$, and $\lambda_2=-\sqrt{\frac{r}{n-r}}$. We claim that $M=M_{r,n-r}$ for some $1\leq r\leq n$. In fact, for any $\alpha>n+1$, from~\eqref{okj}, we have
\[
0=d h_{ij}^\alpha=h_{ij}^{n+1}\omega_{\alpha,n+1}
\]
for any $i,j$. Thus $\omega_{\alpha,n+1}=0$. 
As a consequence, we have 
\[
d\omega_{\alpha\beta}=-\sum_{\gamma>n+1}\omega_{\alpha\gamma}\wedge\omega_{\gamma\beta}
\]
for any $\alpha,\beta>n+1$. Thus locally we can change the frame of the normal bundle such that $\omega_{\alpha\beta}\equiv 0$ for $\alpha,\beta>n+1$.  Evidently, $M$ has to be in some of the totally geodesic submanifold $S^{n+1}$. By using ~\cite{chern-d-k}*{page 68}, we conclude that $M=M_{r,n-r}$ for some $r$.

{\bf Case 3.}
If $A^3=\cdots=A^m=0$ and $A^1, A^2$ are 
\begin{equation}
A_1=\lambda\begin{pmatrix}
\frac{1}{\sqrt 2} &0\\0&-\frac{1}{\sqrt 2}\\
&&0\\
&&&\ddots\\
&&&&0
\end{pmatrix},\quad
A_2=\mu\begin{pmatrix}
0&\frac{1}{\sqrt 2}\\\frac{1}{\sqrt 2}&0\\
&&0\\
&&&\ddots\\
&&&&0
\end{pmatrix},
\end{equation}
$\lambda\mu\neq 0$, and
\[
h^{\alpha}_{ijk}=0
\]
for $\alpha=n+1,n+2$. From ~\eqref{okj}, we have
\[
\frac{1}{\sqrt 2}d\lambda=d h_{11}^{n+1}=0.
\]
Thus $\lambda$ is a constant. Similarly, by computing $dh_{12}^{n+2}$, we know that $\mu$ is also a constant. Without loss of generality, we assume that $\lambda^2\geq\mu^2$.
Let $n>2,j\geq 3$. Then from $0=d h_{1j}^{n+1}=h_{11}^{n+1}\omega_{1j}$ we conclude  $\omega_{1j}=0$ if $j\geq 3$. Similarly, $\omega_{2j}=0$ for $j\geq 3$. Thus by the structure equations, we have
\[
0=d\omega_{1j}=\omega_1\wedge\omega_j,
\]
a contradiction if $n>2$. Thus $n=2$. By computing $dh_{12}^{n+1}$ and $d h_{11}^{n+2}$ using~\eqref{okj}, we get $2\lambda\omega_{12}+\mu\omega_{n+1,n+2}=0$ and
 $2\mu\omega_{12}+\lambda\omega_{n+1,n+2}=0$. Thus $\lambda^2=\mu^2=2/3$. Since $||\sigma||^2=4/3$,  by~\cites{lianmin,xusenlin}, we conclude that $M$ is a Veronese surface. 
 
 {\bf Case 4.} We assume that $n\geq 3$ and $\lambda\mu\neq 0$. Otherwise, we are back to Case 2 or Case 3. We will prove that $M$ doesn't exist. Using (2) of Lemma~\ref{lem2}, we have
\begin{align}
&\omega_{n+1,1}=(n-1)\lambda\omega_1; \label{df-a}\\
&\omega_{n+1,j}=-\lambda\omega_j,\quad 2\leq j\leq n;\label{df-b}\\
&\omega_{\alpha,1}=\mu\omega_{\alpha-n},\quad\alpha>n+1;\label{df-c}\\
&\omega_{\alpha,j}=\delta_{j,\alpha-n}\mu\omega_1,\quad \alpha>n+1, j\geq 2.\label{df-d}
\end{align}
Furthermore, we have
\[
h^{n+1}_{ijk}=0.
\]
$\lambda,\mu$ are presumably local functions, however, since
\[
d\lambda=dh_{11}^{n+1}=0,  
\]
$\lambda$ must be a constant. On the other hand, by
\begin{equation}\label{df-1}
n(n-1)\lambda^2+2n\mu^2=||\sigma||^2+\lambda_2\equiv n,
\end{equation}
 $\mu$ is also a constant.

By differentiating ~\eqref{df-b} using the structure equations, we have
\begin{equation}\label{1001}
-\omega_{n+1,1}\wedge\omega_{1j}-\sum_{k\geq 2}\omega_{n+1,k}\wedge\omega_{kj}-\sum_{\alpha>n+1}\omega_{n+1,\alpha}\wedge\omega_{\alpha,j}=-\lambda d\omega_j.
\end{equation}
However, by~\eqref{df-b}, we have
\[
-\sum_{k\geq 2}\omega_{n+1,k}\wedge\omega_{kj}=\sum_{k\geq 2}\lambda\omega_{jk}\wedge\omega_k
=-\lambda d\omega_j.
\]
 Thus from~\eqref{1001} we conclude
 \[
 \mu\omega_{n+1,n+j}-(n-1)\lambda\omega_{1j}=a_j\omega_1
 \]
 for local smooth functions $a_j$. Let $j\geq 2$, from
 \begin{equation}\label{1002}
 0=dh_{1j}^{n+1}=n\lambda\omega_{1j}-\mu\omega_{n+1,n+j},
 \end{equation}
 we conclude that
 \begin{equation}\label{1003}
 \omega_{1j}=b_j\omega_1,
 \end{equation}
 where $b_j=a_j/\lambda$. Let $j\neq\alpha-n$ and $j\geq 2, \alpha>n+1$. Then from ~\eqref{df-d}, we have
 \[
 0=d\omega_{\alpha j}=-\omega_{\alpha1 }\wedge\omega_{1j}-\sum_{k\geq 2}\omega_{\alpha k}\wedge \omega_{kj}-\omega_{\alpha,n+1}\wedge\omega_{n+1,j}
 -\sum_{\beta>n+1}\omega_{\alpha\beta}\wedge\omega_{\beta j}.
 \]
 
  Thus for $k\neq j$, $k,j\geq 2$, using ~\eqref{df-c},~\eqref{df-d},~\eqref{1002},~\eqref{1003},we have
  \[
  \mu b_j\omega_k\wedge\omega_1-\frac{n\lambda^2}{\mu} b_k\omega_j\wedge\omega_1
  +\mu\omega_1\wedge(\omega_{k,j}-\omega_{n+k,n+j})=0.
  \]
  The third term of the above equation is skew-symmetric with respect to $k,j$. Thus we have
  \begin{equation}\label{ikl}
 \left (\mu-\frac{n\lambda^2}{\mu}\right) (b_k\omega_j\wedge\omega_1+b_j\omega_k\wedge\omega_1)=0.
  \end{equation}
  If 
$
  \mu-\frac{n\lambda^2}{\mu}=0
$,
using ~\eqref{df-c},~\eqref{df-d}, and ~\eqref{1002},  we have
\[
\mu d\omega_1=d\omega_{n+j,j}=0.
\]
Since $d\omega_1=-\omega_{1j}\wedge\omega_j=-b_j\omega_1\wedge\omega_j$,
we have
$\omega_{1j}=0$. If $\mu-\frac{n\lambda^2}{\mu}\neq 0$, then we also have $\omega_{1j}=0$ by~\eqref{ikl}. Using the structure equations, we have
    \[
  0=d\omega_{1j}=\frac 12 R_{1jkl}\omega_s\wedge\omega_t=(1-(n-1)\lambda^2-\mu^2)\omega_1\wedge\omega_j.
  \]
  Thus
  \[
  1-(n-1)\lambda^2-\mu^2=0.
  \]
  Combining with ~\eqref{df-1} we get $\mu=0$, a contradiction.
  
  The theorem is proved.
  
  \qed
  
\begin{cor}
Let $\xi$ be a unit normal vector of $M$ and let $A^\xi$ be the second fundamental form in the $\xi$ direction. If $m\geq 2$ and
\[
0< ||\sigma||^2+\max\,||A^\xi||^2\leq n,
\]
then $M$ has to be the  Veronese surface.
\end{cor}

\qed

 The quantity $||\sigma||^2+\lambda_2$ might be the right object to study pinching theorems. To justify this, we  end this section by making the following conjecture:

\begin{conj}\label{pt}
Let $M$ be an $n$-dimensional minimal submanifold in $S^{n+m}$. If $||\sigma||^2+\lambda_2$ is a constant and if
\[
||\sigma||^2+\lambda_2>n,
\]
then there is a constant $\eps(n,m)>0$ such that
\[
||\sigma||^2+\lambda_2>n+\eps(n,m).
\]
\end{conj}

If $m=1$, this conjecture was proved in~\cite{peng-terng}.

\section{Proof of  the B\"ottcher-Wenzel Conjecture}

In this section, we prove Conjecture~\ref{bw-conj}.

 We fix $X$ and assume that $||X||=1$. Let $V={\mathfrak g}{\mathfrak l}(n,\mathbb R)$. Define a linear map 
\[
T: V\to V,\quad Y\mapsto [X^T,[X,Y]],
\]
where $X^T$ is the transpose of $X$.
Then we have

\begin{lemma}\label{lem1}
$T$ is a semi-positive definite symmetric linear transformation of $V$.
\end{lemma}

{\bf Proof.} This is a straightforward computation
\[
\langle Y_1,[X^T,[X,Y_2]]\rangle=\langle[X,Y_1],[X,Y_2]\rangle
=\langle [X^T,[X,Y_1]],Y_2\rangle.
\]

Obviously $T$ is semi-positive.

\qed

The conjecture is equivalent to the statement that the maximum eigenvalue of $T$ is not more than  $2$.

Let $\alpha$ be the maximum eigenvalue of $T$. Then $\alpha>0$. Let $Y$ be an eigenvector of $T$ with respect to $\alpha$. Then we have
\[
T(Y)=\alpha Y.
\]
A straightforward computation gives
\[
T([X^T,Y^T])=\alpha [X^T,Y^T].
\]

We claim that $Y$ and $Y_1=[X^T,Y^T]$ are linearly  independent: in fact, we have $Y_1\neq 0$, and $\langle Y,Y_1\rangle=0$. Obviously $[X^T,Y^T]\neq 0$. Thus, we have the  following conclusion:

\begin{prop}
\label{prop2}
The multiplicity of  eigenvalue $\alpha$ is at least $2$.
\end{prop}
\qed

Let~\footnote{The method of singular decomposition was first used in~\cite{bw}.} 
\[
X=Q_1\Lambda Q_2
\]
be the singular decomposition of $X$, where $Q_1,Q_2$ are orthogonal matrices and $\Lambda$ is a diagonal matrix. Let
\[
B=Q_2YQ_2^{-1},\quad C=Q_1^{-1}YQ_1.
\]
Then we have
\[
||[X,Y]||^2=||\Lambda B-C\Lambda||^2.
\]
Let
\[
\Lambda=\begin{pmatrix}
s_1\\&\ddots\\&&s_n
\end{pmatrix}.
\]
Without loss of generality, we assume that $s_1^2\geq\cdots \geq s_n^2$. Since $||X||=1$, we have
\[
s_1^2+\cdots+s_n^2=1.
\]

Assume that $s_1^2\leq 1/2$. Then we have
\begin{equation}\label{bw-4}
||\Lambda B-C\Lambda||^2=\sum_{i,j=1}^n(s_ib_{ij}-s_j c_{ij})^2
\leq \sum_{i,j=1}^n 2(b_{ij}^2+c_{ij}^2)s_1^2\leq 2||Y||^2,
\end{equation}
where $B=(b_{ij})$ and $C=(c_{ij})$.
Thus in this case, the conjecture is  true.
Now assume that $s_1^2>1/2$. By Proposition~\ref{prop2}, we can find an eigenvector $Y$ of $T$ such that: 1). $||Y||=1$;  2). the corresponding $b_{11}=0$; and 3). $T(Y)=\alpha Y$.

Conjecture~\ref{bw-conj} follows from the inequality:
\begin{equation}\label{43}
||[X,Y]||^2\leq 2,
\end{equation}
where $Y$ is the particular eigenvector chosen above.

Since  $b_{11}=0$, we have
\[
||\Lambda B-C\Lambda||^2=c_{11}^2 s_1^2
+\sum_{i=2}^n(s_ib_{i1}-s_1c_{i1})^2
+\sum_{j=2}^n (s_1b_{1j}-s_jc_{1j})^2+\Delta_1,
\]
where
\[
\Delta_1=\sum_{i,j=2}^n(s_ib_{ij}-s_j c_{ij})^2.
\]
Since $s_2^2\leq 1/2$, we have
\[
\Delta_1\leq \sum_{i,j=2}^n (b_{ij}^2+c_{ij}^2).
\]
Thus~\eqref{43} is implied by
the following inequality
\begin{equation}\label{ser}
c_{11}^2s_1^2
+\sum_{i=2}^n(s_ib_{i1}-s_1c_{i1})^2
+\sum_{j=2}^n (s_1b_{1j}-s_jc_{1j})^2
\leq \Delta+\sum_{i=2}^n b_{i1}^2+\sum_{j=2}^n c_{1j}^2,
\end{equation}
where 
\[
\Delta=\sum_{i=2}^n b_{1i}^2+\sum_{j=2}^n c_{j1}^2+c_{11}^2.
\]

We consider the matrix
\[
P=
\begin{pmatrix}
\Delta& -b_{12}c_{12}-b_{21}c_{21}&\cdots &-b_{1n}c_{1n}-b_{n1}c_{n1}\\
-b_{12}c_{12}-b_{21}c_{21}& b_{21}^2+c_{12}^2\\
\vdots&&\ddots\\
-b_{1n}c_{1n}-b_{n1}c_{n1}&&&b_{n1}^2+c_{1n}^2
\end{pmatrix}.
\]
Inequality ~\eqref{ser}  is equivalent to that
 the maximum eigenvalue of the above matrix is no more than $\Delta+\sum_{i=2}^n b_{i1}^2+\sum_{j=2}^n c_{1j}^2$. To prove the fact, we let
\[
y=\Delta+\sum_{i=2}^n b_{i1}^2+\sum_{j=2}^n c_{1j}^2+\eps
\]
for $\eps>0$. We have
\[
\det(yI-P)=\prod_{i=2}^n (y-b_{i1}^2-c_{1i}^2)\left(
y-\Delta-\sum_{i=2}^n \frac{(b_{1i}c_{1i}+b_{i1}c_{i1})^2}{y-b_{i1}^2-c_{1i}^2}\right).
\]
Let
\[
\beta=\underset{i>1}{{\rm Max}}\, (b_{i1}^2+c_{1i}^2).
\]
Then we have
\[
y-\Delta-\sum_{i=2}^n \frac{(b_{1i}c_{1i}+b_{i1}c_{i1})^2}{y-b_{i1}^2-c_{1i}^2}\geq \beta+\eps-\beta\sum_{i=2}^n\frac{b_{1i}^2+c_{i1}^2}{\sum_{i=2}^n( b_{1i}^2+c_{i1}^2)}>0.
\]
The conjecture is proved.

\qed

\end{document}